\documentclass[11pt]{amsart}
\begin{document}
%%%%%%%%%%%%%%%%%%%%%%%%%%%%%%%
%%% ALIAS
%%%%%%%%%%%%%%%%%%%%%%%%%%%%%%%

\newtheorem{proposition}{Proposition}[section]
\newtheorem{theorem}{Theorem}[section]
\newtheorem{lemma}{Lemma}[section]
\newtheorem{corollary}{Corollary}[section]
\newtheorem{definition}{Definition}[section]
\newtheorem{conjecture}{Conjecture}[section]

\numberwithin{equation}{section}
\renewcommand{\theequation}{\thesection.\arabic{equation}}

%%%%%%%%%%%%%%%%%%%%%%%%%%%%%%%%%% new commands %%%%%%%%%%%%%%%%%%%%%%
\newcommand{\bp}{\mathfrak{b}}
\newcommand{\xib}{\boldsymbol{\xi}}
\newcommand{\me}{\mathrm{e}} 
\newcommand{\dif}{\mathrm{d}}
\newcommand{\rdm}{\mathfrak{m}}
\newcommand{\jactr}{\nabla\! {\mathbf X}^{\dagger}}
\newcommand{\jac}{\nabla\! {\mathbf X}}
\newcommand{\dersigmatr}{\nabla\!\mbox{\boldmath$\sigma$}_j^{\dagger}}
\newcommand{\dersigma}{\nabla\!\mbox{\boldmath$\sigma$}_j}
\newcommand{\derdrift}{\nabla\!{\mathbf m}}
\newcommand{\derdriftstrat}{\nabla\!{\mathbf q}}
%%%%%%%%%%%%%%%%%%%%%%%%%%%%%%%%%%%%%%%%%%%%%%%%%%%%%%%%%%%%%%%%%%%%%%%%
\newcommand{\forfilt}[2]{\ensuremath{\mathbb{B}_{#1}^{#2}}}  % forward filtration
\newcommand{\backfilt}[2]{\ensuremath{\mathbb{B}_{#1}^{#2}}} % backward filtration
\newcommand{\filt}[1]{\forfilt{}{#1}}
\newcommand{\bdown}{\ensuremath{b^{\downarrow}_T}}            % backward wiener process
\newcommand{\bback}{B}
\newcommand{\lifetimesharp}{\ensuremath{\mathfrak{e}}}        % life time of process sharp
\newcommand{\invmeas}{\ensuremath{\dif\Pi}}                   % Pi invariant measure of X up  
\newcommand{\xinf}{\ensuremath{\mathfrak{X}}}                 % X_\infty
\newcommand{\real}[1]{\ensuremath{\mathbb{R}^{#1}}}           % field of real numbers
\newcommand{\gen}{\ensuremath{\mathcal{G}}}                   % generator symbol
\newcommand{\trace}{\mathtt{tr}}                              % trace
\newcommand{\strat}{\!\bullet\!}
\newcommand{\firstder}{\ensuremath{\frac{\dif}{\dif x}}}             % 1st derivative operator
\newcommand{\secondder}{\ensuremath{\frac{\dif^2}{\dif x^2}}}        % 2nd derivative operator
\newcommand{\genup}{\ensuremath{\frac {\sigma^2}{2}\frac{\dif^2}{\dif x^2}+\left(m+\frac 1{2}\sigma\sigma'\right)\frac{\dif}{\dif x}}}    % generator of the original diffusion
\newcommand{\gensharp}{\ensuremath{\frac {\sigma^2}{2}\frac{\dif^2}{\dif x^2}-\left(m-\frac 1{2}\sigma\sigma'\right)\frac{\dif}{\dif x}}}       % generator of sharp process
%%%%%%%%%%%%%%%%%%%%%%%%%% space functions %%%%%%%%%%%%%%%%%%%%%%%%%%%%%%%%%%%%%%
\newcommand{\Czero}{\ensuremath{\mathcal{C}}}       % Continuous functions
\newcommand{\Czerob}{\ensuremath{\mathcal{C}_b}}    % Continuous bounded functions
\newcommand{\Ctwo}{\ensuremath{\mathcal{C}^2}}    % Continuous functions and differentiable up to 2
\newcommand{\Ctwob}{\ensuremath{\mathcal{C}^2_b}} % Continuous bounded functions and differentiable up to 2
\newcommand{\smooth}{\ensuremath{\mathcal{C}^\infty}} % smooth functions
\newcommand{\Cvan}{\ensuremath{\mathcal{C}_\infty}} % continuos functions vanishing at infinity
%%%%%%%%%%%%%%%%%%%%%%%%%%%%%%%%%%%%%%%%% processes %%%%%%%%%%%%%%%%%%%%%%%%%%%%
\newcommand{\proup}[1]{\ensuremath{X_{#1}}}                % diffusion up
\newcommand{\proupb}[1]{\ensuremath{\mathbf{X}_{#1}}}      % diffusion up boldface
\newcommand{\prodown}[1]{\ensuremath{X^{\downarrow}_{#1}}} % X down  
\newcommand{\prodownb}[1]{\ensuremath{\mathbf{X}^{\downarrow}_{#1}}} % X down boldface  
\newcommand{\prointer}{\ensuremath{X^T}}                  
\newcommand{\prointermediate}[1]{\ensuremath{X^{T}_{#1}}}    % X^T intermediate
\newcommand{\prosharp}[1]{\ensuremath{X^{\sharp}_{#1}}}        % X^# sharp 
\newcommand{\prosharpb}[1]{\ensuremath{\mathbf{X}^{\sharp}_{#1}}}   % X^# sharp boldface 
\newcommand{\prosharpinv}[1]{\ensuremath{X^{\sharp -1}_{#1}}}   % X^# sharp inverse 
\newcommand{\proback}[1]{\ensuremath{\mathfrak{z}_{#1}}}              % X pullback
%%%%%%%%%%%%%%%%%%%%%%%%%%%%%%%%%%%%%%%%%% flows %%%%%%%%%%%%%%%%%%%%%%%%%%%%%%%
\newcommand{\flowup}[1]{\ensuremath{X_{#1}}}                   % X up 
\newcommand{\flowdown}[1]{\ensuremath{X^{\downarrow}_{#1}}}    % X down  
\newcommand{\flowintermediate}[1]{\ensuremath{X^{T}_{#1}}}     % X^T intermediate
\newcommand{\flowsharp}[1]{\ensuremath{X^{\sharp}_{#1}}}       % X^# sharp 
\newcommand{\flowsharpinv}[1]{\ensuremath{X^{\sharp -1}_{#1}}} % X^# sharp inverse 
\newcommand{\invmap}[1]{\ensuremath{\phi_{#1}}}

%%%%%%%%%%%%  TITLE  %%%%%%%%%%%%%%
\title[On the invariant measure]{On the invariant measure of a positive recurrent diffusion in \real{}}
\author[M. Baldini]{Michele Baldini$^\dag$}
\thanks{$^\dag$Department of Mathematics, Northwestern University,
2033 Sheridan Road, Evanston, IL 60208. Email: \texttt{mlb@math.northwestern.edu}}
\thanks{This paper was written at New York University, as part of a Ph.D. thesis. The author is extremely grateful to his advisor Prof. Henry McKean for his invaluable guidance and contagious enthusiasm.}
%%%%%%%%%%%  ABSTRACT  %%%%%%%%%%%%%
\begin{abstract}
Given an one-dimensional positive recurrent diffusion governed by the Stratonovich SDE
\[
X_t=x+\int_0^t\sigma(X_s)\strat\dif b(s)+\int_0^t m(X_s)\dif s,
\]
we show that the associated stochastic flow of diffeomorphisms focuses as fast as $
\mathrm{exp}\left(-2t\int_{\real{}}\frac{m^2}{\sigma^2}\dif\Pi\right)$,
where $\invmeas$ is the finite stationary measure. Moreover, if the drift is reversed and the diffeomorphism is inverted, then the path function so produced tends, independently of its starting point, to a single (random) point whose distribution is $\invmeas$. 
Applications to stationary solutions of \proup{t}, asymptotic behavior of solutions of SPDEs and random attractors are offered.
\end{abstract}
\maketitle

%%%%%%%%%%%%%%%%%%%%%%%%%%%%%%%%%%%%%%%
%%%%% INTRODUCTION
%%%%%%%%%%%%%%%%%% %%%%%%%%%%%%%%%%%%%%%
\section{Introduction}
\label{introduction}
A well-known result \cite{hasminskii},\cite{pinsky} states that a diffusion in $\real{d}$ possesses an unique invariant density if and only if it is positive recurrent, i.e. if it visits almost surely a  neighborhood of every point $z$ in a finite time $\tau_z$ and ${\bf E}[\tau_z]<\infty$. Heuristically, positive recurrence on \real{d} is synonymous of a strong restoring drift that makes the diffusion  spend most of the time near the origin, thus producing an invariant measure of finite total mass. In general, the invariant density is identified as the positive solution of $\gen^{*}u=0$, where $\gen^{*}$ is the formal adjoint of the infinitesimal generator, but this might be a daunting task. In this manuscript we offer an alternative interpretation of such stationary state in the case $d=1$, with the hope that the theory can be extended to more advanced systems, in which the the computation of the invariant measure is still an open problem (for instance, the stochastic Burgers equation). We announce that partial results are already obtained in the case $d\ge 2$, but this is left for future publication.

Consider a positive recurrent diffusion \proup{t} in \real{}, solution of the SDE  
\begin{equation}
X_t=x+\int_0^t\sigma(X_s)\strat\dif b(s)+\int_0^t m(X_s)\dif s,
\end{equation}
and let \prosharp{t} be the diffusion governed by the same $\sigma$ but opposite drift, i.e.
\begin{equation}
\prosharp{t}=x+\int_0^t\sigma(\prosharp{s})\strat\dif b(s)-\int_0^t m(\prosharp{s})\dif s.
\end{equation}
Here $\sigma$ and $m$ are smooth real functions and $\strat$ denotes the Stratonovich integral. For almost every Brownian path and for all $t\ge 0$, it is folklore that the maps $x\mapsto\proup{t}(x)$ and $x\mapsto\prosharp{t}(x)$ may be regarded as local diffeomorphisms of \real{} \cite{kunita},\cite{protter}. These two flows, and in particular their asymptotic behavior, will be the main objects of our investigation. 

The paper is organized as follows. In section \ref{preliminaries} we will engage in  a study of the ergodic properties of \prosharp{t} and its associated stochastic flow. 
In section \ref{focusing} we will prove that the distance (in the natural scale) between two solutions of \proup{t} driven by the same Brownian motion decreases exponentially fast to zero; i.e. let $s$ be the scale function, then for every $x,y\in\real{}$, 
\begin{equation}\label{rate}
{\bf P}\left[\lim_{t\to\infty}\frac{1}{t}\ln\big|s\big(\proup{t}(x)\big)-s\big(\proup{t}(y)\big)\big|=-2\int_{\real{}}\frac{m^2}{\sigma^2}\dif\Pi\right]=1,
\end{equation}
where $\invmeas$ is the invariant measure of \proup{t}. Interestingly, the factor $-2\int_{\real{}}\frac{m^2}{\sigma^2}\dif\Pi$ provides also an upper bound for the spectral gap of \proup{t}.
In section \ref{constructionxdown} and \ref{stagnationpoint} we will show that there exists a random variable $\xinf$ (depending upon the whole history of the Brownian path) whose law can be identified with the invariant measure $\invmeas$. Such random variable is constructed as the limit of a new path function $\prodown{t}$ which converges almost surely as $t\uparrow\infty$.
The process $\prodown{t}$ admits a representation in terms of 
a stochastic flow of diffeomorphisms: its inverse map is a transient diffusion whose drift is exactly that of the original diffusion with the sign changed. It turns out that $\xinf$ satisfies:
\begin{itemize}
\item ${\bf P}\left[\liminf_{t\to\infty}|\prosharp{t}(\xinf)|<\infty\right]=1$.
\item ${\bf P}\left[\lim_{t\to\infty}|\prosharp{t}(x)|=\infty\text{ for all }x\neq\xinf\right]=1$.
\item ${\bf P}\left[\lim_{t\to\infty}\prodown{t}(x)=\xinf,\text{ for all }x\right]=1$.
\end{itemize}
In section \ref{ornstein} we offer an example, while
sections \ref{stationary_solutions} through \ref{random_attractors} are devoted to applications. We will prove that 
\begin{itemize}
\item $\xinf$ can be used to construct stationary solutions of both \proup{t} and \prosharp{t}.
\item For every smooth function $f$, the process $f(\prodown{t})$ with initial condition $f(\prodown{0})=f(x)$ is the unique solution of 
\begin{equation}
f(\prodown{t})=f(x)-\sigma(x)\int_0^t\frac{\partial f(\prodown{s})}{\partial x}\dif b(s)+\int_0^t\gen f(\prodown{s})\,\dif s,
\end{equation}
where $\gen f(\prodown{t})$ represents the infinitesimal generator of $\proup{t}$ applied to the whole function $f\big(\prodown{t}\big)$.
\item $\xinf$ can be identified with the random attractor of the corresponding random dynamical system.
\end{itemize}
{\bf Notation}: The symbol $\bp$ denotes a generic Brownian path. We will work with a standard Wiener space $(\Omega,\forfilt{}{},{\bf P})$, enlarged so that $\real{}\times\Omega\ni(t,\mathfrak{b})\mapsto b(t)$ is a two-sided Brownian motion, i.e. two independent standard Brownian motion $b_1(t)_{t\ge 0}$ and $b_2(t)_{t\ge 0}$ connected back to back at $b(0)=0$ such that $b(t)=b_1(t)$ if $t\ge 0$ and $b(t)=b_2(-t)$ if $t\le 0$. The shift $\theta_t\bp=\theta_t b(s):=b(s+t)-b(t), t\in\real{}$ is a measure preserving and ergodic operator on the space of Brownian paths. $\int f\dif b$ is an It\^{o} integral and $\int f\strat \dif b$ is a Stratonovich integral, while $\int f\overleftarrow{\dif b}$
and $\int f\strat\overleftarrow{\dif b}$ denote the corresponding backward integrals (see \cite{kunita} for the definitions). $\mathcal{C}^\alpha\ (\mathcal{C}_b^\alpha)$ is the space of (bounded) functions $\alpha$-times differentiable with continuity. In particular, $\mathcal{C}^\infty$ is called the set of {\it smooth} functions. 
%%%%%%%%%%%%%%%%%%%%%%%%%%%%%%%%%%%%%%%%%
%%% PRELIMINARIES
%%%%%%%%%%%%%%%%%%%%%%%%%%%%%%%%%%%%%%%%%
\section{Preliminaries}
\label{preliminaries}
\subsection{The main characters}
Let $\sigma:\real{}\mapsto (0,\infty)$ and $m:\real{}\mapsto\real{}$ be smooth functions such that
\begin{equation}\label{positive_rec_condition1}
\int_{-\infty}^{0}\frac{1}{\sigma(x)}\mathrm{exp}\left(-\int_{0}^{x}\frac{2m}{\sigma^2}\right)\dif x=\int_{0}^{\infty}\frac{1}{\sigma(x)}\mathrm{exp}\left(-\int_{0}^{x}\frac{2m}{\sigma^2}\right)\dif x=\infty,
\end{equation}
\begin{equation}\label{positive_rec_condition2}
\Lambda\stackrel{\mathrm{def}}{=}\int_{\real{}}\frac{1}{\sigma(x)}\mathrm{exp}\left(\int_{0}^x \frac{2m}{\sigma^2}\right)\dif x<\infty.
\end{equation}
For such regular $\sigma$ and $m$ and for every $x\in\real{}$, it easy to prove \cite{mckean_stoch_int} that there exists
\begin{itemize}
\item a diffusion \proup{t} defined up to an explosion time $\mathfrak{e}(x,\mathfrak{b})$ that satisfies the Stratonovich SDE
\begin{equation}\label{recurrent_diffusion}
X_t=x+\int_0^t\sigma(X_s)\strat \dif b(s)+\int_0^tm(X_s)\dif s,
\end{equation}
\item a diffusion \prosharp{t} (read {\it sharp}) defined up to an explosion time $\mathfrak{e}^\sharp(x,\mathfrak{b})$ that satisfies the Stratonovich SDE
\begin{equation}\label{sharp_transient_diffusion}
X^{\sharp}_t=x+\int_0^t\sigma(X^{\sharp}_s)\strat \dif b(s)-\int_0^tm(X^{\sharp}_s)\dif s.
\end{equation}
\end{itemize}
The infinitesimal generator, the transition probability function, the semigroup and the invariant measure associated to $X_t$ are denoted respectively $\gen$, $\,p_t(x,\cdot)$, $\mathrm{exp}(t\gen)$ and $\invmeas$, whereas the counterparts for the {\it sharp} diffusion are indicated with a superscript $\sharp$. Let us introduce the positive function
\begin{equation}
\psi^2(x)\doteq\Lambda\,\sigma(x)\mathrm{exp}\left(-\int_0^x \frac{2m}{\sigma^2}\right)
\end{equation}
and the measures $\dif\mu(x)= 2\psi^{-2}(x)\dif x$, $\dif s(x)=\psi^2(x)\sigma^{-2}(x)\dif x$.  A straighforward computation produces
\begin{eqnarray}
& & \gen\equiv\genup=\frac 1{2}\psi^2\firstder\left(\frac{\sigma^2}{\psi^2}\firstder\right)=\frac{\dif}{\dif\mu}\frac{\dif}{\dif s},\\[3mm]
& & \gen^{\sharp}\equiv\gensharp= \frac{\sigma^2}{\psi^2}\firstder\left(\frac{\psi^2}{2}\firstder\right)=\frac{\dif}{\dif s}\frac{\dif}{\dif \mu},
\end{eqnarray}
showing that the speed measures and the scale measures of \proup{t} and \prosharp{t} are ``duals'' of each other. Hence, $\psi^{-2}$ is the unique invariant (probability) density of $\proup{t}$ and similarly $\psi^2\sigma^{-2}$ is the unique invariant density (of infinite total mass) of $\prosharp{t}$. 

In the next theorem, we study the boundary and the ergodic property of these two processes.  \begin{theorem} 
Let \proup{t} and \prosharp{t} be as above. Then 
\begin{enumerate}\renewcommand{\theenumi}{\alph{enumi}} 
\item \proup{t} is a positive recurrent diffusion.
\item \prosharp{t} is a transient diffusion, whose boundary points $+\infty$ and $-\infty$ are {\it non-entrance}. In particular,
\begin{eqnarray}
& & {\bf P}_x\left[\lim_{t\to\infty}\prosharp{t}=+\infty\right]=\int_{-\infty}^x \invmeas,\\
& & {\bf P}_x\left[\lim_{t\to\infty}\prosharp{t}=-\infty\right]=\int_x^\infty \invmeas.
\end{eqnarray}
\end{enumerate} 
\end{theorem}
\begin{proof}
(a) Conditions \eqref{positive_rec_condition1} and \eqref{positive_rec_condition2} are necessary and sufficient to guarantee positive recurrence of \proup{t} \cite{pinsky}. (b) It is known \cite{pinsky} that the property of transience is equivalent to the existence, for every $x\in\real{}$, of a finite Green measure $G^\sharp(x,\dif y)=\int_0^\infty p^\sharp_t(x,\dif y)\dif t$ which can be identified as the limit of $G_n^\sharp(x,\dif y)$, where
\begin{equation}
G_n^\sharp(x,\dif y)=
\begin{cases}
2\frac{\displaystyle s^\sharp(x)-s^\sharp(-n)}{\displaystyle s^\sharp(n)-s^\sharp(-n)}\left(s^\sharp(n)-s^\sharp(y)\right)\dif\mu^\sharp(y) & \text{when }y\ge x,\\[4mm]
2\frac{\displaystyle s^\sharp(n)-s^\sharp(x)}{\displaystyle s^\sharp(n)-s^\sharp(-n)}\left(s^\sharp(y)-s^\sharp(-n)\right)\dif\mu^\sharp(y) & \text{when }y\le x.
\end{cases}
\end{equation}
Hence, \prosharp{t} is transient if and only if $\lim_{n\to\infty}s^\sharp(n)$ and $\lim_{n\to\infty}s^\sharp(-n)$ exist and are finite. This follows immediately from $\int_{-\infty}^{\infty}\dif s^\sharp=\int_{-\infty}^{\infty}\dif \mu=2$.
Let $z\in\real{}$ and let $\tau^\sharp_z$ be the first hitting time of $z$. Then, from the formulas for exit probabilities, 
\begin{equation}
{\bf P}_x\left[\lim_{t\to\infty}\prosharp{t}=+\infty\right]={\bf P}_x\Big[\tau_{\infty}<\tau_{-\infty}\Big]=\frac{\int_{-\infty}^x \dif s^\sharp}{\int_{-\infty}^\infty \dif s^\sharp}=\int_{-\infty}^x \invmeas.
\end{equation}
To complete the proof, it remains to study the boundary. According to Feller's classification \cite{feller}, $+\infty$ is a non-entrance point if $\int_0^\infty \dif s^\sharp(y)\int_y^\infty\dif\mu^\sharp(z)=\infty$.
Since $s(+\infty)=\infty$ and integrating by parts 
\begin{equation}
\int_0^\infty \dif s^\sharp(y)\int_y^\infty\dif\mu^\sharp(z)=\int_0^\infty d\mu(y)\int_y^\infty \dif s(z)=\infty.
\end{equation}
An analogous result holds for $-\infty$.
\end{proof}
{\bf Comment:} Because the non-entrance property, $+\infty$ and $-\infty$ are the cimitery points for \prosharp{t}. After its death, $\prosharp{t}$ is buried at infinity for good.

%%%%%%%%%%%%%%%%%%%%%%%%%%%%%%%%%%%%%%%%%%
%% Flows of diffeomorphisms 
%%%%%%%%%%%%%%%%%%%%%%%%%%%%%%%%%%%%%%%%%%
\subsection{Stochastic flows of diffeomorphisms} 
In virtue of the recurrence property, both $+\infty$ and $-\infty$ are {\it non-exit} boundary points for $\proup{t}$. Feller \cite{feller} showed that this implies {\it conservativeness}, i.e. for each $x$, $\proup{t}$ has infinite life-time with probability 1. Moreover, in one dimension it is a simple computation to prove that conservativeness implies {\it strictly conservativeness}, i.e. $\proup{t}$ has infinite life-time for all $x$ with probability 1. According to Kunita \cite{kunita}, this is sufficient to guarantee that, after a null set in $\Omega$ has been weeded out, the family of maps $\{\proup{t}(x):t\ge 0\}$ defines a stochastic flow of local diffeomorphisms, i.e. for any Brownian path $\bp$
\begin{enumerate}\renewcommand{\theenumi}{\alph{enumi}}
\item $\real{}\ni x\mapsto\proup{t}(x,\mathfrak{b})$ is a {\bf into} diffeomorphism for all $t\ge0$ and the inverse map is continuous in $t$ and smooth in $x$.
\item $\proup{t+\tau}(x,\bp)=\proup{t}(\cdot,\theta_{\tau}\bp)\circ\proup{\tau}(x,\bp)$ for all non-negative $t,\tau$.
\end{enumerate}

Define $\mathfrak{R}_t(\mathfrak{b})=\inf\{x:\flowsharp{t}(x,\bp)=+\infty\}$, 
$\mathfrak{L}_t(\mathfrak{b})=\sup\{x:\flowsharp{t}(x,\bp)=-\infty\}$, 
and let $D_t(\mathfrak{b})\stackrel{def}{=}\{x: t<\mathfrak{e}^\sharp(x,\mathfrak{b})\}=(\mathfrak{L}_t(\mathfrak{b}), \mathfrak{R}_t(\mathfrak{b}))$. Clearly, for any $t\le\tau$ we have $D_t(\mathfrak{b})\subseteq D_{\tau}(\mathfrak{b})$. In  virtue of the transience property, both $+\infty$ and $-\infty$ are {\it non-entrance} boundary points for \prosharp{t}.  According to \cite{kunita}, this is sufficient to guarantee that, after a null set in $\Omega$ has been weeded out, the family of maps $\{\prosharp{t}(x):t\ge 0\}$ defines a stochastic flow of local diffeomorphisms, i.e. for any Brownian path $\bp$
\begin{enumerate}\renewcommand{\theenumi}{\alph{enumi}}
\item $D_t(\mathfrak{b})\ni x\mapsto\prosharp{t}(x,\mathfrak{b})$ is a {\bf onto} diffeomorphism for all $t\ge0$ and the inverse map is continuous in $t$ and smooth in $x$.
\item $\prosharp{t+\tau}(x,\bp)=\prosharp{t}(\cdot,\theta_{\tau}\bp)\circ\prosharp{\tau}(x,\bp)$ holds on $D_{t+\tau}(\mathfrak{b})$ for all non-negative $t,\tau$.
\end{enumerate}

In what follows, we will use the symbols \proup{t} and \prosharp{t} to indicate the stochastic processes and $\proup{t}(x,\bp)$ and $\prosharp{t}(x,\bp)$ to indicate the associated flows of local diffeomorphisms. Often, we will omit $\bp$. We will work in the "normalized" space of paths $\mathcal{W}\subseteq\Omega$,
such that both $\proup{t}(x,\mathfrak{b})$ and $\prosharp{t}(x,\mathfrak{b})$ are well defined for all $\bp\in\mathcal{W}$. Clearly, ${\bf P}(\mathcal{W})=1$. 
\subsection{Stagnation and invariant points}
From the definition of transience we have ${\bf P}_x[\lim_{t\to\infty}|\prosharp{t}|=\infty]=1$ for every $x\in\real{}$.
However, for each Brownian path, there might be some $x$ such that $\prosharp{t}(x,\bp)$ does not run out to infinity.
\begin{definition} For every Brownian path define  
\begin{equation}
\mathcal{S}(\mathfrak{b})=\{x:\liminf_{t\to\infty}|\flowsharp{t}(x,\bp)|<\infty\}.
\end{equation}
$\mathcal{S}(\mathfrak{b})$ is called the {\bf stagnant} set. 
A random variable $Y(\mathfrak{b})$ is called a {\bf stagnation point} if 
$Y(\mathfrak{b})\in\mathcal{S}(\mathfrak{b})$ almost surely. 
\end{definition}
\begin{theorem} For almost every Brownian path, the stagnant set  $\mathcal{S}(\mathfrak{b})$ contains at most one point.
\end{theorem}
\begin{proof} For every $x$ define $K(x)=\{\bp:\liminf_{t\to\infty}|\flowsharp{t}(x,\bp)|<\infty\}$ and let ${\bf 1}_{K(x)}(\bp)$ be its indicator function. Then ${\bf E}[{\bf 1}_{K(x)}(\bp)]=0$ for all $x$. Integrating on $x$ and using Fubini's theorem, we obtain ${\bf 1}_{K(x)}(\bp)=0$ for almost every $\bp$ and a.e. in $x$, concluding that that the set $\mathcal{S}(\bp)$ has zero Lebesgue measure almost surely. By virtue of monotonicity of the map $\prosharp{t}(x)$, if there were two distinct points $x_1<x_2$ such that
$|\flowsharp{t}(x_1)|\not\rightarrow\infty$ and $|\flowsharp{t}(x_2)|\not\rightarrow\infty$ as 
$t\uparrow\infty$, then all the trajectories originating from the interval
$[x_1,x_2]$ would be ``sandwiched'' between $\flowsharp{t}(x_1)$ and $\flowsharp{t}(x_2)$. Hence the whole interval $[x_1,x_2]$ would belong to $\mathcal{S}(\bp)$, contradicting the fact that its measure is zero.  
\end{proof}

Now introduce the random variables $R(\mathfrak{b})=\inf\{x:\lim_{t\to\infty}\flowsharp{t}(x,\bp)=+\infty\}$, $L(\mathfrak{b})=\sup\{x:\lim_{t\to\infty}\flowsharp{t}(x,\bp)=-\infty\}$, and observe 
$L(\bp)=R(\bp)$ almost surely. In fact, by definition $L(\bp)\le R(\bp)$. If the inequality were strict, $|\flowsharp{t}(x)|\not\rightarrow\infty$ as $t\uparrow\infty$ for any $L(\bp)<x<R(\bp)$, contradicting the transience of $\prosharp{t}$. Define $\xinf(\bp)=L(\bp)(=R(\bp))$. At this point the reader might wonder whether $\xinf$ is a true stagnation point: the question will be answered positively in section \ref{stagnationpoint}. To summarize what we have so far
\begin{theorem} 
\label{xinfexistence}
There exists a unique random variable $\xinf\in\forfilt{0}{\infty}$ such that, almost surely, 
\begin{eqnarray}
& &\lim_{t\to\infty}\flowsharp{t}(x,\bp)=+\infty\quad\mbox{for all }x>\xinf,\nonumber\\[1mm]
& &\lim_{t\to\infty}\flowsharp{t}(x,\bp)=-\infty\quad\mbox{for all }x<\xinf.
\end{eqnarray}
\end{theorem} 
\begin{definition} 
A random variable $Y(\mathfrak{b})$ is called an {\bf invariant point} if
\begin{equation}
{\bf P}\left[\flowsharp{t}\big(\cdot,\bp\big)\circ Y(\mathfrak{b})=Y\big(\theta_t\mathfrak{b}\big)\text{ for all }t\ge 0\right]=1.
\end{equation}
\end{definition}
\begin{lemma}\label{invariant_implies_stagnation}
If $Y(\bp)$ is an invariant point, then is also a stagnation point.
\end{lemma}
\begin{proof} Let $Y(\bp)$ be an invariant point, then a moment of reflection shows that the process $\prosharp{t}(Y)\equiv\prosharp{t}(Y(\bp),\bp)$ is stationary end ergodic. Let $\Lambda$ be the set of Brownian paths such that  $\liminf_{t\to\infty}|\flowsharp{t}(Y)|=\infty$.  Let $K$ be any compact subset of \real{} and let ${\bf 1}_K(\cdot)$ be its indicator function. Then for each Brownian path in  $\Lambda$, we have
\begin{equation}
\lim_{T\to\infty}\frac{1}{T}\int_0^T {\bf 1}_K\big(\flowsharp{t}(Y)\big)\dif t=0.
\end{equation}
On the other hand, combining Birkhoff's Theorem with ergodicity, we can always adjust the set $K$ such that almost surely  
\begin{equation}
\lim_{T\to\infty}\frac{1}{T}\int_0^T {\bf 1}_K\big(\prosharp{t}(Y)\big)\dif t={\bf E}\Big[{\bf 1}_K(Y(\bp))\Big]=\int_K{\bf P}(Y(\bp)\in\dif x)>0.
\end{equation}
Therefore ${\bf P}(\Lambda)=0$.
\end{proof}

%%%%%%%%%%%%%%%%%%%%%%%%%%%%%%%%%%%%%%%%%%%%
%%% FOCUSING X UP
%%%%%%%%%%%%%%%%%%%%%%%%%%%%%%%%%%%%%%%%%%%%
\section{Asymptotic focusing of the flow $\proup{t}(x)$.}
\label{focusing}
In \cite{hasminskii_focusing}, Has'minskii proved that any solution of a one-dimensional stochastic differential equation, describing a positive recurrent diffusion, is stable with probability one in the metric given by the natural scales. More precisely, he showed that for any $x,y\in\real{}$ 
\begin{equation}
{\bf P}\Big[\lim_{t\to\infty}\Big|s\big(\proup{t}(x)\big)-s\big(\proup{t}(y)\big)\Big|=0\Big]=1,
\end{equation}
where $s(\cdot)$ is the scale function. Below, we offer a different proof of this stability and we improve the result showing that the pathwise rate of convergence is exponentially fast.  
Let us introduce the smooth map $\frac{\partial\proup{t}(x)}{\partial x}$ from $\real{}$ into $\real{}$. The associated process is governed by the SDE
\begin{equation}
\dif \frac{\partial\proup{t}}{\partial x}=\sigma'(X_t)\frac{\partial\proup{t}}{\partial x}\dif b(t)+q'(X_t)\frac{\partial\proup{t}}{\partial x}\dif t,\quad\frac{\partial\proup{0}}{\partial x}=1,
\end{equation}
where $q$ is the modified drift $m+\frac{1}{2}\sigma\sigma'$. 
\begin{lemma}\label{convergence_gradient_proup}
Let \proup{t} be a positive recurrent diffusion as in \eqref{recurrent_diffusion} and let $s(x)=\int_0^x\frac{\psi^2}{\sigma^2}$ be its scale function. If $\gamma\stackrel{def}{=} 2\int_\real{}\frac{m^2}{\sigma^2}\invmeas<\infty$, then for every $x$
\begin{equation}
{\bf P}_x\left[\lim_{t\to\infty}\frac{1}{t}\ln\frac{\partial s\big(\proup{t}\big)}{\partial x}=-\gamma\right]=1.
\end{equation}
\end{lemma}
\begin{proof}
An application of It\^{o}'s lemma gives 
\begin{equation}
\dif \frac{\partial s\big(\proup{t}\big)}{\partial x}=-\frac{2m}{\sigma}\big(\proup{t}\big)\frac{\partial s\big(\proup{t}\big)}{\partial x}\dif b(t), 
\end{equation}
leading to 
\begin{equation}\label{exp_supermart_gradient}
\frac{\partial s\big(\proup{t}\big)}{\partial x}=s'(x)\mathrm{exp}\left(-2\int_0^t \frac{m}{\sigma}\left(\proup{\tau}\right)\dif b(\tau)-2\int_0^t \frac{m^2}{\sigma^2}\left(\proup{\tau}\right)\dif \tau\right).
\end{equation}
For every $x$, the RHS of \eqref{exp_supermart_gradient} is a positive local martingale and by virtue of Fatou's lemma is also a supermartingale. The exponent is in the form of 
$w(S_t)-\frac{1}{2}S_t$ where $w(S_t)$ is a new Brownian motion running with the random clock $
S_t=\int_0^t \frac{4m^2}{\sigma^2}\big(\proup{\tau}\big)\dif \tau$.
Recurrence of \proup{t} guarantees that the clock grows to infinity and an application of the law of the iterated logarithm to the Brownian integral produces for almost all Brownian paths and for any $\epsilon>0$,
\begin{equation}
-(1+\epsilon)\sqrt{2S_t\ln\ln S_t}\le w(S_t)\le (1+\epsilon)\sqrt{2S_t\ln\ln S_t}\quad\text{for }t\uparrow\infty.
\end{equation}
Hence,
\begin{equation}
s'(x)\me^{-(1+\epsilon)\sqrt{2S_t\ln\ln S_t}-\frac{1}{2}S_t}\le\frac{\partial s\big(\proup{t}(x)\big)}{\partial x}\le s'(x)\me^{(1+\epsilon)\sqrt{2S_t\ln\ln S_t}-\frac{1}{2}S_t}\quad\text{for }t\uparrow\infty.
\end{equation}
The proof is completed using the ergodic theorem.
\end{proof}
\begin{theorem}\label{convergence_metric_scale}
Let \proup{t} be a positive recurrent diffusion as in \eqref{recurrent_diffusion} and let $s(x)=\int_0^x\frac{\psi^2}{\sigma^2}$ be its scale function. For any $a<b\in\real{}$ write $A_t\stackrel{def}{=}s\big(\proup{t}(a)\big)$ and $B_t\stackrel{def}{=}s\big(\proup{t}(b)\big)$. If $\gamma<\infty$, then 
\begin{equation}
{\bf P}\Big[\lim_{t\to\infty}\Big|B_t-A_t\Big|=0\Big]=1.
\end{equation}
\end{theorem}
\begin{proof}
Lemma \ref{convergence_gradient_proup} can be used to show that $B_t-A_t$ is a positive supermartingale and therefore there exists a positive random variable $\zeta_\infty$ such that $\lim_{t\to\infty}(B_t-A_t)=\zeta_\infty$ almost surely.
Define $\zeta_t=B_t-A_t$, $f(x)=\arctan(x)$ and using the ergodic theorem observe 
\begin{equation}
0<\frac{1}{T}\int_0^T\Big[f(B_t)-f(A_t)\Big]dt\xrightarrow{T\to\infty}0\quad\text{a.s.}
\end{equation}
An application of the mean-value theorem and bounded convergence, shows that
\begin{equation}
\left|\frac{1}{T}\int_0^T\Big[f(A_t+\zeta_t)-f(A_t+\zeta_\infty)\Big]dt\right|
\end{equation}
vanishes almost surely as $T\to\infty$. Hence, from the trivial identity  
\begin{multline}
\frac{1}{T}\int_0^T\Big[f(B_t)-f(A_t)\Big]dt\\
=\frac{1}{T}\int_0^T\Big[f(A_t+\zeta_t)-f(A_t+\zeta_\infty)\Big]dt+\frac{1}{T}\int_0^T\Big[f(A_t+\zeta_\infty)-f(A_t)\Big]dt,
\end{multline}
we have
\begin{equation}\label{identity_gihman}
{\bf P}\left[\lim_{t\to\infty}\frac{1}{T}\int_0^T\Big[f(A_t+\zeta_\infty)-f(A_t)\Big]dt=0\right]=1.
\end{equation}
It easy to show that the family $F_T(x)$ is equicontinuos and $
\lim_{T\to\infty}F_T(x)=\int_{\real{}}\Big[f(z+x)-f(z)\Big]\invmeas(z)$,
uniformly on each finite interval. Then, 
$\lim_{T\to\infty}F_T(\zeta_\infty)=\int_{\real{}}\Big[f(z+\zeta_\infty)-f(z)\Big]\invmeas(z)$
and a comparison with \eqref{identity_gihman} produces $
\int_{\real{}}\Big[f(z+\zeta_\infty)-f(z)\Big]\invmeas(z)=0$ almost surely.
Because $f(z+\zeta_\infty)-f(z)$ is non-negative and continuous, then $f(z+\zeta_\infty)-f(z)=0$ for all $z$, indicating $\zeta_\infty=0$. This nice trick is adapted from \cite{gihman}.
\end{proof}
%%%%%%%%%%%%%%%%%%%%%%%%%%%%%%%%%%%%%%
%%% pathwise rate of convergence
%%%%%%%%%%%%%%%%%%%%%%%%%%%%%%%%%%%%%%
\begin{theorem}
Under the same assumptions of theorem \ref{convergence_metric_scale},
\begin{equation}
{\bf P}\left[\lim_{t\to\infty}\frac{1}{t}\ln\left(B_t-A_t\right)=-\gamma\right]=1.
\end{equation}
\end{theorem}
\begin{proof}{\bf Step 1: The clock goes up.} Define 
\begin{equation}
F(x,y)=\frac{\frac{\psi^2}{\sigma}(y)-\frac{\psi^2}{\sigma}(x)}{\int_x^y \frac{\psi^2}{\sigma^2}}
\end{equation}
and using It\^{o}'s lemma observe that 
\begin{equation}\label{martingale_scale_difference}
\frac{B_t-A_t}{s(a)-s(b)}=\mathrm{exp}\left(\int_0^t F\Big(\proup{\tau}(a),\proup{\tau}(b)\Big)\dif b(\tau)-\frac{1}{2}\int_0^t F^2\Big(\proup{\tau}(a),\proup{\tau}(b)\Big)\dif\tau\right).
\end{equation}
In virtue of theorem \ref{convergence_metric_scale},
$ \lim_{t\to\infty}\int_0^t F^2\Big(\proup{\tau}(a),\proup{\tau}(b)\Big)\dif\tau=\infty$ almost surely and hence 
$B_t-A_t\approx\mathrm{exp}(-\frac{1}{2}M_t)$ for $t\uparrow\infty$, where $M_t\stackrel{def}{=} \int_0^t F^2\Big(\proup{\tau}(a),\proup{\tau}(b)\Big)\dif\tau$.

{\bf Step 2: Bound from above.} Using the proof of lemma \ref{convergence_gradient_proup} as a model, for almost all Brownian paths, and for any $\epsilon>0$,
\begin{equation}\label{constrains_for_focusing}
\me^{-(1+\epsilon)\sqrt{2M_t\ln\ln M_t}-\frac{1}{2}M_t}\le\frac{B_t-A_t}
{s(b)-s(a)}\le\me^{(1+\epsilon)\sqrt{2M_t\ln\ln M_t}-\frac{1}{2}M_t}\quad\text{for }t\uparrow\infty.
\end{equation}
Define $\rho(x)=(\frac{\psi^2}{\sigma}\circ s^{-1})'(x)$. An application of the mean-value theorem gives
\begin{equation}
\frac{\psi^2}{\sigma}\big(\proup{\tau}(b)\big)-\frac{\psi^2}{\sigma}\big(\proup{\tau}(a)\big)=\rho\big(C_\tau\big)\Big(B_\tau-A_\tau\Big)
\end{equation}
for some $C_\tau\in\big(A_\tau,B_\tau\big)$, leading to 
$M_t=\int_0^t \rho^2\big(C_\tau\big)\dif\tau$.
Pick a subsequence of non-negative functions $F_n$ with bounded first derivative such that $|F'_n|\le L_n$ for a certain constant $L_n$ and $F_n\uparrow\rho^2$. Using the mean-value theorem,
\begin{equation}
\int_0^t F_n(C_\tau)\dif\tau=\int_0^t F_n(A_\tau)\dif\tau+\int_0^t F^{\prime}_n(D_\tau)\left(C_\tau-A_\tau\right)\dif\tau
\end{equation}
for some $D_\tau\in(A_\tau,C_\tau)$. Because $M_t\ge \int_0^t F_n(C_\tau)\dif\tau$, then
\begin{multline}
\mathrm{exp}\left(-\frac{M_t}{2}\right)\le\mathrm{exp}\left(-\frac{1}{2}\int_0^t F_n(A_\tau)\dif\tau+\frac{1}{2}\left|\int_0^t F^{\prime}_n(D_\tau)\left(C_\tau-A_\tau\right)\dif\tau\right|\right)\\[3mm]
\le\mathrm{exp}\left(-\frac{1}{2}t\int_{\real{}}(F_n\circ s)\invmeas-t\times o(1)+\frac{1}{2}L_nt\int_0^1 \left(C_{t\tau}-A_{t\tau}\right)\dif\tau\right).
\end{multline}
 Therefore, $B_t-A_t\le (s(b)-s(a))\me^{-\frac{1}{2}M_t\left(1-o(1)\right)}$, producing
\begin{multline}
\frac{1}{t}\ln\left(B_t-A_t\right)\\
\le\left[-\frac{1}{2}\int_{\real{}}(F_n\circ s)\invmeas-o(1)+\frac{1}{2}L_n\int_0^1 \left(C_{t\tau}-A_{t\tau}\right)\dif\tau\right]\Big(1-o(1)\Big)\\
+\frac{1}{t}\ln\left[s\big(b\big)-s\big(a)\right].
\end{multline}
Now, since $C_{t\tau}-A_{t\tau}\le B_{t\tau}-A_{t\tau}$,  
\begin{equation}
\limsup_{t\to\infty}\frac{1}{t}\ln\left(B_t-A_t\right)\le
-\frac{1}{2}\int_{\real{}}(F_n\circ s)\invmeas.
\end{equation}
In particular,
\begin{equation}
\limsup_{t\to\infty}\frac{1}{t}\ln\left(B_t-A_t\right)\le-\sup_{n\ge 0}\frac{1}{2}\int_{\real{}}(F_n\circ s)\invmeas=-\frac{1}{2}\int_{\real{}}\left(\rho^2\circ s\right)\invmeas=-\gamma.
\end{equation}

{\bf Step 3: Bound from below.} For every $\epsilon>0$, Fatou's lemma gives 
\begin{equation}
\liminf_{t\to\infty}\me^{(\gamma+\epsilon)t}\left(B_t-A_t\right)\ge\int_a^b \liminf_{t\to\infty}\me^{(\gamma+\epsilon)t}\ \frac{\partial s\big(\proup{t}(x)\big)}{\partial x}dx=\infty,
\end{equation}
which permits to claim $\me^{-(\gamma+\epsilon)t}\le\left(B_t-A_t\right)$ for $t\uparrow\infty$. Therefore
$-\gamma-\epsilon\le\liminf_{t\to\infty}\frac{1}{t}\ln\left(B_t-A_t\right)$.

{\bf Step 4: Chain of inequality.} Combining the bounds in step 2 and 3, we obtain for every $\epsilon>0$,
\begin{equation}
-\gamma-\epsilon\le\liminf_{t\to\infty}\frac{1}{t}\ln\left(B_t-A_t\right)\le\limsup_{t\to\infty}\frac{1}{t}\ln\left(B_t-A_t\right)\le-\gamma,
\end{equation}
and the proof is completed letting $\epsilon\downarrow 0$.
\end{proof}
%%%%%%%%%%%%%%%%%%%%%%%%%%%%%%%%%%%%%%%%%%%%%%%%%%%
%%% upper bound for spectral gap of mathcal{G}
%%%%%%%%%%%%%%%%%%%%%%%%%%%%%%%%%%%%%%%%%%%%%%%%%%%
Below we show that the factor $\gamma=2\int_\real{}\frac{m^2}{\sigma^2}\invmeas$ has an interesting interpretation in terms of the spectral gap of the generator $\gen$. Let us recall that the spectral gap $\mathfrak{G}$ of a selfadjoint operator with the eigenvalue corresponding to the ground state equal to zero, is the modulus of the first (non-trivial) eigenvalue. 
\begin{theorem}\label{spectral_gap_less_than_gamma}
Let \proup{t} be a positive recurrent diffusion as in \eqref{recurrent_diffusion}. Define $V(x)=\int_0^x \frac{1}{\sigma(y)}\dif y$ and $\gamma=2\int_\real{}\frac{m^2}{\sigma^2}\invmeas$. If the following three conditions are satisfied,
\begin{enumerate}\renewcommand{\theenumi}{\alph{enumi}}
\item $\int_{\real{}}V^2\invmeas <\infty,$
\item $\int_{\real{}}\sigma^2\invmeas <\infty,$
\item $\gamma<\infty$,
\end{enumerate}
then $\mathfrak{G}\le\gamma$.
\end{theorem}
\begin{proof} By definition, the spectral gap can be characterized as
\begin{equation}
\mathfrak{G}=\inf_{f\in Q}\left\{\frac{1}{2}\int_{\real{}}\frac{(f')^2}{s'}\right\}
\end{equation}
where
$Q=\left\{f\in L^2(\real{},\invmeas)\text{ and smooth, }\int_{\real{}}f\invmeas=0,\ \int_{\real{}}f^2\invmeas=1\right\}$.
Since $\gamma$ can be rewritten as
$
\gamma=\frac{1}{2}\int_\real{}\frac{1}{s'}\left(\frac{\sigma'}{\sigma}-2\frac{\psi'}{\psi}\right)^2$,
it suffices to prove
\begin{equation}
\inf_{f\in Q}\left\{\frac{1}{2}\int_{\real{}}\frac{(f')^2}{s'}\right\}\le \frac{1}{2}\int_\real{}\frac{1}{s'}\left(\frac{\sigma'}{\sigma}-2\frac{\psi'}{\psi}\right)^2.
\end{equation}
Define $F(x)=\alpha(V(x)+\beta)$ where 
$\beta=-\int_\real{}V\invmeas$ and 
$\frac{1}{\alpha^2}=\int_\real{}\big(V+\beta\big)^2\invmeas$.
The definition has been made to guarantee $F\in Q$; also $\int_\real{}(V+\beta)'\frac{\sigma}{\psi^2}=1$ by choice of $V$. In virtue of the assumed conditions, an application of the Cauchy-Schwarz inequality shows that the integral 
\begin{equation}
\int_\real{}(V+\beta)\left(\frac{\sigma}{\psi^2}\right)'=\int_\real{}\frac{V+\beta}{\psi}\left(\frac{2m}{\sigma\psi}\right)
\end{equation}
exists and is finite. Therefore, it is permissible to integrate by parts in\\ $\int_\real{}(V+\beta)'\frac{\sigma}{\psi^2}$ obtaining
\begin{eqnarray}
1=\left[\int_\real{}(V+\beta)'\frac{\sigma}{\psi^2}\right]^2 & = & \left[\int_\real{}(V+\beta)\frac{\sigma}{\psi^2}\left(\ln\frac{\sigma}{\psi^2}\right)^{'}\right]^2\label{inequality_for_the_gap}\nonumber\\[2mm]
& \le & \int_\real{}(V+\beta)^2\invmeas\times\int_\real{}\frac{1}{s'}\left(\frac{\sigma'}{\sigma}-2\frac{\psi'}{\psi}\right)^2=\frac{2\gamma}{\alpha^2}\nonumber.
\end{eqnarray}

Finally, 
\begin{equation}
\mathfrak{G}\le\frac{1}{2}\int_{\real{}}\frac{(F')^2}{s'}=\frac{\alpha^2}{2}\times\left[\int_\real{}(V+\beta)'\frac{\sigma}{\psi^2}\right]^2\le \frac{1}{2}\int_\real{}\frac{1}{s'}\left(\frac{\sigma'}{\sigma}-2\frac{\psi'}{\psi}\right)^2=\gamma.
\end{equation}
\end{proof}
{\bf Remark:} The reader might wonder if there exists a choice of $\sigma$ and $m$ such that $\mathfrak{G}=\gamma$. If $\sigma$ is constant, it is a simple exercise to prove that the equality holds if and only if the drift is affine.
%%%%%%%%%%%%%%%%%%%%%%%%%%%%%%%%%%%%%%%%%%%%%
%%% CONSTRUCTION X DOWN
%%%%%%%%%%%%%%%%%%%%%%%%%%%%%%%%%%%%%%%%%%%%%
\section{The process \prodown{t}: construction and asymptotic behavior}
\label{constructionxdown}
For any $T>0$, consider the Brownian motion $\bdown(t)_{0\le t\le T}$ defined as 
$\bdown(t)= b(T-t)-b(T)$. Let $\{\prointermediate{t}: 0\le t\le T\}$ be the solution of the SDE  
\begin{equation}
\dif X_t^T=\sigma(X_t^T)\strat \dif\bdown(t) +m(X_t^T)\dif t,\quad\prointermediate{0}=x. 
\end{equation}
and define $\prodown{T}=\prointermediate{T}$. The family $\{\prodown{t}: t\ge 0\}$ describes a new (non-Markovian) process starting at $x$ (read {\it down}). Note that for each $t\ge 0$ separately, the process $\prodown{t}$ is equal in law to \proup{t}. 
\begin{lemma}\label{inverse_diffeomorphism}
Let $\prosharp{t}(x,\bp)$ be the flow associated to
\[
X^{\sharp}_t=x+\int_0^t\sigma(X^{\sharp}_s)\strat \dif b(s)-\int_0^tm(X^{\sharp}_s)\dif s.
\]
Then for every $\bp$, the inverse map satisfies $
\prosharpinv{t}(x,\bp)=\prodown{t}(x,\bp)$ for all $t\ge 0$.
\end{lemma}
\begin{proof} For every $\bp$ and all $0\le t\le T$, the composition of maps rule gives  
$\prosharpinv{T}\big(x,\bp\big)=\prosharpinv{t}\big(\cdot,\bp\big)\circ\prosharpinv{T-t}\big(x,\theta_t\bp\big)$, which leads to $\prosharp{t}\big(\cdot,\bp\big)\circ\prosharpinv{T}\big(x,\bp\big)=\prosharpinv{T-t}\big(x,\theta_t\bp\big)$. 
Denoting $\prosharpinv{T}\big(x,\bp\big)=y$, a straightforward but tedious computation shows
\begin{equation}
\int_0^T \sigma\left(\flowsharp{t}\big(y,\bp\big)\right)\strat \dif b(t)
=-\int_0^T \sigma\left(\prosharpinv{t}\big(x,\theta_{T-t}\bp\big)\right)\strat \dif\bdown(t),
\end{equation}
and
\begin{equation}
-\int_0^T m\left(\flowsharp{t}\big(y,\bp\big)\right)\dif t=-\int_0^T m\left(\flowsharpinv{t}\big(x,\theta_{T-t}\bp\big)\right)\dif t.
\end{equation}
Combining the above identities, we obtain
\begin{equation}
\prosharpinv{T}(x,\bp)=x+\int_0^T\sigma\left(\prosharpinv{t}(x,\theta_{T-t}\bp)\right)\strat \dif\bdown(t)+\int_0^T m\left(\prosharpinv{t}(x,\theta_{T-t}\bp)\right)\dif t
\end{equation}
and a comparison with the flow $X^T_t(x)$, produces $\prosharpinv{T}(x,\bp)=\prointermediate{T}(x,\bp)\equiv\prodown{T}
(x,\bp)$.
\end{proof}

%%%%%%%%%%%%%%%%%%%%%%%%%%%%%%%%%%%%%%%%%%%%%
%%% ASYMPTOTIC BEHAVIOR B DOWN
%%%%%%%%%%%%%%%%%%%%%%%%%%%%%%%%%%%%%%%%%%%%%
Now we finally present the main theorem. Recall from section \ref{preliminaries} that 
\begin{equation}
\xinf=\inf\{x:\lim_{t\to\infty}\flowsharp{t}(x)=\infty\}=\sup\{x:\lim_{t\to\infty}\flowsharp{t}(x)=-\infty\}.
\end{equation}
\begin{theorem}\label{convergenceflowdown} 
The flow $\flowdown{t}(x)$ converges to $\xinf$ almost surely, i.e.
\begin{equation}
{\bf P}\left[\lim_{t\to\infty}\prodown{t}(x)=\xinf,\text{ for all }x\right]=1.
\end{equation}
\end{theorem}
\begin{proof} 
From standard results, either (a) ${\bf P}\big(\mathfrak{e}^\sharp(x,\bp)=\infty\big)=1$ for all $x$ or (b) ${\bf P}\big(\mathfrak{e}^\sharp(x,\bp)=\infty\big)<1$ for all $x$ \cite{mckean_stoch_int}.\\  
%%%%%%%%%%%%%%%%%%%%%%%%%%%%%%%%%%%%%%
%%%% PROOF: NON EXPLOSIVE CASE
%%%%%%%%%%%%%%%%%%%%%%%%%%%%%%%%%%%%%%
(a) Because $+\infty$ and $-\infty$ are {\it non-exit} boundary points for \prosharp{t}, for every Brownian path $\flowsharp{t}(x,\bp)$ and the inverse map 
$\flowdown{t}(x,\bp)$ are flows of global diffeomorphisms of \real{}, for all $t\ge 0$. Now, for any $\epsilon>0$, define $\mathfrak{l}_\epsilon=\xinf -\epsilon$ 
and $\mathfrak{r}_\epsilon=\xinf +\epsilon$. In virtue of theorem \ref{xinfexistence}, for any $x\in\real{}$ 
we have $\flowsharp{t}(\mathfrak{l}_\epsilon)\le x\le\flowsharp{t}(\mathfrak{r}_\epsilon)$ for $t\uparrow\infty$.
Therefore $\mathfrak{l}_\epsilon=\flowdown{t}(\cdot)\circ\flowsharp{t}(\mathfrak{l}_\epsilon)\le \flowdown{t}(x)\le \flowdown{t}(\cdot)\circ\flowsharp{t}
(\mathfrak{r}_\epsilon)=\mathfrak{r}_\epsilon$
holds for sufficiently large time, and thus
\[ 
\mathfrak{l}_\epsilon\le\liminf_{t\to\infty} \flowdown{t}(x)\le\limsup_{t\to\infty}
\flowdown{t}(x)\le \mathfrak{r}_\epsilon. 
\] 
The proof is completed letting $\epsilon\downarrow 0$.\\
%%%%%%%%%%%%%%%%%%%%%%%%%%%%%%%%%%%%%%%%%%%%%%
%%%  PROOF: EXPLOSIVE CASE
%%%%%%%%%%%%%%%%%%%%%%%%%%%%%%%%%%%%%%%%%%%%%%
(b) For every Brownian path, $\flowsharp{t}(x,\bp):(\mathfrak{L}_t,\mathfrak{R}_t)\mapsto\real{}$ is a onto diffeomorphism for all $t\ge 0$ with inverse map $\prodown{t}(x,\bp)$ (see section \ref{preliminaries} for the definition of $\mathfrak{L}_t$ and $\mathfrak{R}_t$). Clearly $\mathfrak{L}_t\le\xinf\le\mathfrak{R}_t$ and since $\mathfrak{e}^\sharp(\mathfrak{L}_t,\bp)=\mathfrak{e}^\sharp(\mathfrak{R}_t,\bp)=t$, then $\mathfrak{L}_t$ is strictly increasing and $\mathfrak{R}_t$ is strictly decreasing as $t\uparrow\infty$. Therefore both $\mathfrak{L}_t$ and $\mathfrak{R}_t$ converge to $\xinf$ as $t\uparrow\infty$.
Because $\mathfrak{L}_t\le\flowdown{t}(x)\le\mathfrak{R}_t$ holds almost surely for all $x$, then the proof is completed taking $t\uparrow\infty$.
\end{proof}
\begin{corollary}\label{corollary_xdown_convergence_AS}
The probability distribution of $\xinf$ is the invariant measure $\invmeas$.
\end{corollary}
\begin{proof}
Let $f\in\Czerob(\real{})$. For every $t\ge 0$ separately, \prodown{t} and \proup{t} are equally distributed. Hence, 
for any $x\in\real{}$
\begin{equation}\label{xinfdistribution} 
{\bf E}_x\left [ f\big(\prodown{t}\big) \right ]= {\bf E}_x\Big[ f\big(\proup{t}\big)
\Big]= \left(\me^{t\mathcal{G}}f\right)(x). 
\end{equation}
The LHS of \eqref{xinfdistribution} converges to ${\bf E}\left[f(\xinf)\right]$ as $t\uparrow\infty$, meanwhile the ergodic theorem for a Feller semigroup guarantees $\left(\me^{t\mathcal{G}}f\right)(x)\xrightarrow{t\to\infty} \int_{\real{}} f\invmeas$. Therefore $
{\bf E}\Big[ f(\xinf) \Big]=\int_{\real{}}\! f\invmeas$ for any arbitrary (continuous, bounded) function $f$.
\end{proof}
%%%%%%%%%%%%%%%%%%%%%%%%%%%%%%%%%%%%%%
%%%% Xinf IS THE STAGNATION POINT
%%%%%%%%%%%%%%%%%%%%%%%%%%%%%%%%%%%%%%
\section{$\xinf$ is the stagnation point}
\label{stagnationpoint}
In section \ref{preliminaries}, we showed that for almost every Brownian paths, the stagnation set  $\mathcal{S}(\bp)=\{x:\liminf_{t\to\infty}|\flowsharp{t}(x,\bp)|<\infty\}$
is either empty or contains one point. In a second step, we introduced the (random) point $\xinf$ and we established that if $\mathcal{S}(\bp)$ is not empty, then it must contain such a point. Below, we prove that $\xinf$ is an invariant point, which in virtue of lemma \ref{invariant_implies_stagnation} guarantees the stagnation property.

%%%%%%%%%%%%%%%%%%%% stationarity of X sharp at xinf %%%%%%%%%%%
\begin{theorem}\label{x_sharp_is_stationary}
$\xinf$ is an invariant point.
\end{theorem}\begin{proof}
Because $\flowdown{\tau}\big(x,\bp\big)=\flowdown{t}\big(\cdot,\bp\big)\circ\flowdown{\tau-t}\big(x,\theta_t\bp\big)$, it follows $\flowsharp{t}\big(\cdot,\bp\big)\circ\flowdown{\tau}\big(x, \bp\big)=\flowdown{\tau-t}\big(x, \theta_t\bp\big)$ for every Brownian path and all $0\le t\le\tau$. To complete the proof, let $\tau\uparrow\infty$ and use continuity of the diffeomorphism to show
\begin{equation}\label{prosharp_at_xinf_is_stationary}
\flowsharp{t}\big(\cdot,\bp\big)\circ \xinf(\mathfrak{b})=\xinf\big(\theta_t\mathfrak{b}\big)\text{ for all }t\ge 0,\text{ a.s.}
\end{equation}
\end{proof}
%%%%%%%%%%%%%%%%%%%%%%%%%%%%%%%%%%%%%%%%%%%%%%%%%%%%%%%%%%%%%%%%%%%%%%%%%%%%
%%% Example: Ornstein-Uhlenbeck
%%%%%%%%%%%%%%%%%%%%%%%%%%%%%%%%%%%%%%%%%%%%%%%%%%%%%%%%%%%%%%%%%%%%%%%%%%%%
\section{Example: the Ornstein-Uhlenbeck process}
\label{ornstein}
The O-U process is positive recurrent and satisfies the equation
\begin{equation}\label{OUprocess}
\dif X_t=\dif b(t)-\beta X_t\dif t
\end{equation}
with $\beta>0$. Its solution, with initial condition $x$, reads
\begin{equation}
\proup{t}(x)=x\me^{-\beta t}+\int_0^t \me^{\beta(s-t)}\dif b(s).
\end{equation}
Following the construction in section \ref{constructionxdown},
\begin{equation}
\prodown{t}(x)=x\me^{-\beta t}+\int_0^t \me^{\beta(s-t)}\dif b^{\downarrow}_t(s)=x\me^{-\beta t}-\int_0^t \me^{-\beta \tau}\dif b(\tau).
\end{equation}
Therefore $\xinf=\lim_{t\to\infty}\prodown{t}(x)=-\int_0^\infty \me^{-\beta \tau}\dif b(\tau)$. Compare its law 
with the well-known invariant probability density $\sqrt{(\beta/\pi)}\me^{-\beta x^2}$. Now observe that 
\begin{equation}
\prosharp{t}(\cdot,\bp)\circ\xinf(\bp)=-\me^{\beta t}\int_t^\infty \me^{-\beta s}\dif b(s)
\end{equation}
is Gaussian, continuous and it has the same correlation as the stationary process $-\me^{\beta t}b\left(\me^{-2\beta t}/2\beta\right)$, {\it ergo} is stationary as the doctor ordered.

%%%%% application to stationary solution

\section{Stationary solution of \proup{t} and \prosharp{t}}\label{stationary_solutions}
In section \ref{stagnationpoint}, we used $\xinf$ to prove the existence and uniqueness of a (non-trivial) stationary solution of \prosharp{t}. If the direction of time is reversed, we can also use $\xinf$ to construct a (non-trivial) stationary solution of \proup{t}. 
Consider the SDE \eqref{sharp_transient_diffusion} running backward in the negative semi-axis. Since the adaptedness of the solution must run from right to left, the extension can be made with the introduction of the Stratonovich backward stochastic integral 
\begin{equation}
-\int_{-t}^0 \sigma(\prosharp{s})\strat\overleftarrow{\dif b}(s)=-\int_{-t}^0 \sigma(\prosharp{s}) \overleftarrow{\dif b}(s)+\frac{1}{2}\int_{-t}^0 \sigma(\prosharp{s})\sigma'(\prosharp{s})\dif s
\end{equation}
where $\overleftarrow{\dif b}$ stands for backward It\^{o} integral \cite{kunita}.
Hence, the process $\prosharp{-t}:t\ge 0$ is solution of the SDE
\begin{equation}\label{prosharp_running_backward}
x=\prosharp{-t}-\int_{-t}^0 \sigma(\prosharp{s})\strat\overleftarrow{\dif b}(s)-\int_{-t}^0 m(\prosharp{s})\dif s,
\end{equation}
which, through a simple manipulation, can be rewritten as the forward SDE
\begin{equation}
\prosharp{-t}=x+\int_0^t \sigma(\prosharp{-s})\strat\dif b(-s)+\int_0^t m(\prosharp{-s})\dif s,
\end{equation}
in which the adaptedness of the solution runs from left to right.

Denoting the "rotated" Brownian path with the symbol $\bp^+$, i.e. $\real{}\times\Omega\ni (t,\bp^+)\mapsto b(-t)$,
then the (striking) moral of the story is the following: {\it $\prosharp{t}: (t,\bp)\in\real{-}\times\Omega$ running backward is a version of $\proup{t}: (t,\bp^+)\in\real{+}\times\Omega$ running forward}.

Now, a moment of reflection shows that the whole construction described in sections \ref{constructionxdown} and \ref{stagnationpoint} can be repeated {\it verbatim} in the region $t\le 0$ (picturesquely, imagine to flip the time-line and switch the Brownian motion $b_1(t)$ with $b_2(t)$). Therefore the identity \eqref{prosharp_at_xinf_is_stationary} can be extended to all $t\in\real{}$ almost surely, with the understanding that $(t,\bp)\mapsto b_1(t)$ if $t\ge 0$ and  $(t,\bp)\mapsto b_2(-t)$ if $t\le 0$. This permits to draw the following theorem.

\begin{theorem}
Let \proup{t} and \prosharp{t} be two diffusions as in \eqref{recurrent_diffusion} and \eqref{sharp_transient_diffusion}. Then 
\begin{enumerate}\renewcommand{\theenumi}{\alph{enumi}} \item $\prosharp{-t}(\cdot,\bp^+)\circ\xinf(\bp^+)\equiv\xinf(\theta_{-t}\bp^+)$ is the unique stationary solution
of \proup{t} for $t\ge0$,
\item $\prosharp{t}(\cdot,\bp)\circ\xinf(\bp)\equiv\xinf(\theta_{t}\bp)$ is the unique stationary solution
of \prosharp{t} for $t\ge 0$.
\end{enumerate}
\end{theorem}

%%%%%%%%%%%%% CONNECTIONS WITH SPDEs and PULL_BACK ATTRACTOR
\section{Applications to stochastic partial differential equations (SPDEs)}\label{application_SPDEs}
In the next lemma we introduce a second-order linear SPDE that describes the evolution of the {\it down} process. We focus mainly on the derivation of the equation and we refer to \cite{rozovskii} for general questions of existence and uniqueness.
\begin{lemma}\label{Xdown_satisfies_SPDE}
For any $f\in\Ctwo(\real{})$, the process $f\big(\prodown{t}\big)$ with initial condition $f(x)$ is the unique solution of 
\begin{equation}\label{SPDE}
f\big(\prodown{t}\big)=f(x)-\sigma(x)\int_0^t\frac{\partial f\big(\prodown{s}\big)}{\partial x}\,\dif b(s)+\int_0^t\mathcal{G}f\big(\prodown{s}\big)\,\dif s
\end{equation}
\end{lemma}
\begin{proof} Let $q=m+\frac{1}{2}\sigma\sigma'$ be the modified drift. For all $t\le T$, use the SDEs for the first and second derivative of the map $X_t^T(x)$ to check that $\sigma(x)\frac{\partial\flowintermediate{t}}{\partial x}$ and $\gen\flowintermediate{t}$ satisfy respectively
\begin{equation}
\sigma(x)\frac{\partial\flowintermediate{t}}{\partial x}=\sigma(x)
+\int_0^t\sigma'\big(\prointer_s\big)\sigma(x)\frac{\partial\flowintermediate{s}}{\partial x}\dif\bdown(s)+\int_0^t q'\big(\prointer_s\big)\sigma(x)\frac{\partial\flowintermediate{s}}{\partial x}\dif s\label{equationforcoefficientofdb}
\end{equation}
and
\begin{multline}
\gen\flowintermediate{t}=m(x)+\frac{1}{2}\sigma(x)\sigma'(x)\\
+\int_0^t\sigma'\big(\prointer_s\big)\gen\flowintermediate{s}\dif\bdown(s)+\frac{\sigma^2(x)}{2}\int_0^t\sigma''\big(\prointer_s\big)\left(\frac{\partial\flowintermediate{s}}{\partial x}\right)^2\dif\bdown(s)\label{equationforcoefficientofdt}\\
+\int_0^t q'\big(\prointer_s\big)\gen\flowintermediate{s}\dif s+\frac{\sigma^2(x)}{2}\int_0^t q''\big(\prointer_s\big)\left(\frac{\partial\flowintermediate{s}}{\partial x}\right)^2\dif s.
\end{multline}

Now evaluate $\Delta_h(t,T)\doteq X^{T+h}_{t+h}-\prointermediate{t}$ when $0<h\ll 1$ and $t\le T$, as in
\begin{multline}
\Delta_h(t,T)=\int_0^{h}\sigma\big(X_s^{T+h}\big)\strat\dif b^{\downarrow}_{T+h}(s)+\int_h^{t+h}\sigma\big(X_s^{T+h}\big)\dif b^{\downarrow}_{T+h}(s)\\
+\int_0^{h} m\big(X^{T+h}_s\big)\dif s+\int_h^{t+h} q\big(X^{T+h}_s\big)\dif s\label{expansion}
-\int_0^t\sigma\big(\prointer_s\big)\dif\bdown(s)-\int_0^t q\big(\prointer_s\big)\dif s.
\end{multline}
A tedious manipulation on the stochastic integral shows that the first piece in the RHS of \eqref{expansion} is
\begin{equation}
\sigma(x)b^{\downarrow}_{T+h}(h)+\frac{1}{2}\sigma(x)\sigma'(x)b^{\downarrow 2}_{T+h}(h)+o(h)+o(b^{\downarrow 2}_{T+h}(h)).
\end{equation}
The third term reduces simply to $m(x)h+o(h)$. A Taylor expansion on the second and fifth terms produces
\begin{multline}
\int_0^{t}\left[\sigma\big(X_{s+h}^{T+h}\big)-\sigma\big(\prointer_s\big)\right]\dif\bdown(s)
=\int_0^t\sigma'\big(\prointer_s\big)\Delta_h(s,T)\dif\bdown(s)+\\
\frac{1}{2}\int_0^t\sigma''\big(\prointer_s\big)\Delta_h^2(s,T)\dif\bdown(s)
+\frac{1}{6}\int_0^t\sigma'''(z_2)\Delta_h^3(s,T)\dif\bdown(s), 
\end{multline}
and similarly on the fourth and sixth terms  
\begin{multline}
\int_0^{t}\left[q\big(X^{T+h}_{s+h}\big)-q\big(X^{T}_s\big)\right]\dif s
=\int_0^t q'\big(\prointer_s\big)\Delta_h(s,T)\dif s+\\
\frac{1}{2}\int_0^t q''\big(\prointer_s\big)\Delta_h^2(s,T)\dif s
+\frac{1}{6}\int_0^t q'''(z_3)\Delta_h^3(s,T)\dif s,
\end{multline}
where $z_2,z_3$ are points belonging to the segment $\left[\prointer_s,X^{T+h}_{s+h}\right]$.

Let $h\downarrow 0$ and using the rules of stochastic calculus, we have
\begin{multline}
\Delta_{h}(t,T)=-\sigma(x)\dif b(T)+\left[m(x)+\frac{1}{2}\sigma(x)\sigma'(x)\right]\dif T\\
+\int_0^t\sigma'\big(\prointer_s\big)
\Delta_{h}(s,T)\dif\bdown(s)+\frac{1}{2}\int_0^t\sigma''\big(\prointer_s\big)\Delta_{h}^2(s,T)\dif\bdown(s)\label{eq_for_delta}\\
+\int_0^t q'\big(\prointer_s\big)\Delta_{h}(s,T)\dif s+\frac{1}{2}\int_0^t q''\big(\prointer_s\big)\Delta^2_{h}(s,T)\dif s\\
+\frac{1}{6}\int_0^t\sigma'''(z_2)\Delta_h^3(s,T)\dif\bdown(s)+\frac{1}{6}\int_0^t q'''(z_3)\Delta_h^3(s,T)\dif s. 
\end{multline}
A comparison between terms of the same order on the LHS and the RHS of \eqref{eq_for_delta} shows that $\Delta_{h}(t,T)$ must have both parts of order $\dif b(T)$ and $\dif T$; any other part vanishes effectively.
Therefore, we can write $\Delta_{h}(t,T)$ in the form
\begin{equation}\label{form_of_Delta}
\Delta_{h}(t,T)=\Omega(t,T)\dif b(T)+M(t,T)\dif T,\ \ \ 0\le t\le T.
\end{equation}
Substitute (\ref{form_of_Delta}) in (\ref{eq_for_delta}) and compare the evolution of $\Omega(t,T)$ and $M(t,T)$ with the unique strong solutions of \eqref{equationforcoefficientofdb} and \eqref{equationforcoefficientofdt}. Therefore for all $0\le t\le T$, 
$\Omega(t,T)=-\sigma(x)\frac{\dif\flowintermediate{t}}{\dif x}$ and $M(t,T)=\gen\flowintermediate{t}$ almost surely. Put $t=T$ and observe 
\begin{equation}
\Delta_{\dif T}(T,T)=\dif X^{\downarrow}_T=-\sigma(x)\frac{\partial X^{\downarrow}_T}{\partial x}\dif b(T)+\gen X_T^{\downarrow}\dif T.
\end{equation}
A simple application of It\^{o}'s lemma concludes the proof. 
\end{proof}
\begin{theorem}\label{convergence_of_solution_of_2nd_order_SPDEs}
Let $\sigma:\real{}\mapsto (0,\infty)$ and $\psi:\real{}\mapsto\real{}$ be smooth functions such that 
$\int_{\real{}}\psi^{-2}=1$ and $\int_{-\infty}^0\frac{\psi^2}{\sigma^2}=\int_{0}^\infty\frac{\psi^2}{\sigma^2}=\infty$.
Let $\mathcal{G}$ be the elliptic differential operator $
\mathcal{G}=\frac{1}{2}\psi^2\frac{\dif}{\dif x}\left(\frac{\sigma^2}{\psi^2}\frac{\dif}{\dif x}\right)$.
For any $f\in\Ctwo(\real{})$, then the unique solution of 
\begin{equation}
\mathfrak{u}(t,x)=f(x)-\sigma(x)\int_0^t\frac{\partial \mathfrak{u}}{\partial x}(s,x)\,\dif b(s)+\int_0^t\mathcal{G}\mathfrak{u}(s,x)\,\dif s
\end{equation}
satisfies $
{\bf P}\Big[\lim_{t\to\infty}\mathfrak{u}(t,x)\text{  exists, for all } x\Big]=1$. 
The limit is a finite random variable independent of $x$, whose probability density is $\psi^{-2}$.
\end{theorem}
\begin{proof} The integrability conditions on the coefficients $\sigma$ and $\psi$ guarantee that the elliptic operator $\gen$ is associated to a positive recurrent diffusion. Then the proof follows immediately from 
lemma \ref{Xdown_satisfies_SPDE} and theorem \ref{convergenceflowdown}.
\end{proof}
To the best of our knowledge, the above theorem presents new results on long-time behavior of random fields that are solutions of SPDEs in unbounded domains. Chueshov and Vuillermot \cite{vuillermot} have investigated a similar problem. Their parabolic equation is studied in a bounded domain and has been subjected to a homogeneous white noise whose coefficient depends on the solution itself but not on its gradient.
%%%%%%%%%%%%%%%%%%%%%%%%%%%%%%%%%%%%%%%%%%%%%%%%%%
%%% random attractors
%%%%%%%%%%%%%%%%%%%%%%%%%%%%%%%%%%%%%%%%%%%%%%%%%%
\section{Applications to random attractors}\label{random_attractors}
A well-known results \cite{arnold} states that a Stratonovich stochastic differential equation as \eqref{recurrent_diffusion} generates a two-sided random dynamical system (RDS) with respect to the canonical shift $\theta_t$ on the Wiener space. In connection with the study of attractors and invariant sets, special attention is reserved to the {\it pull-back} convergence of the RDS as opposed to the usual forward convergence. This is to observe the values at time $t=0$ of solutions that originate from the same fixed point, progressively earlier in time. Below, we will define the pull-back convergence at the level of SDE and in a second step, we will prove that it is equivalent to the convergence of the process \prodown{t}, provided that $\bp$ and $\bp^+$ are swapped. 

Pick a $T\ge 0$ and consider the diffusion \proup{t} running forward in the region $-T\le t\le 0$. Since the Brownian motion travels from right to left, the process will be adapted to the filtration $\forfilt{-T}{t}=\sigma(b(u)-b(-T):-T\le u\le t)$ and it is solution of the backward SDE (backward with respect of the negative semi-axis)
\begin{equation}\label{proup_running_forward_t<0}
\proup{t}=x+\int_{-T}^t \sigma(\proup{s})\strat\dif b(s)+\int_{-T}^t m(\proup{s})\dif s.
\end{equation}
Put $t=0$ and define $\proback{T}=\proup{0}$. Now repeat the same operation for all $T\ge 0$, maintaining the same initial condition $\proup{-T}=x$. We say that family $\proback{t}:t\ge 0$ starting at $\proback{0}=x$ is the {\it pull-back} process associated to \proup{t}.
\begin{definition}
If there exists a random variable $Y(\mathfrak{b})\in\forfilt{-\infty}{0}$ such that 
\begin{equation}
{\bf P}[\lim_{t\to\infty}\proback{t}=Y\text{ for all }x]=1
\end{equation}
we say that the  diffusion \proup{t} converges to $Y$ in a pull-back sense.
\end{definition}
It is a simple exercise to prove 
\begin{lemma}\label{spde_for_pull-back_process} 
For any $f\in\Ctwo(\real{})$, the process $f(\proback{t})$ starting at $f(x)$ is the unique solution of 
\begin{equation}\label{SPDE_pullback}\nonumber
f(\proback{t})=f(x)-\sigma(x)\int_0^t\frac{\partial f(\proback{s})}{\partial x}\,\dif b(-s)+\int_0^t\mathcal{G}f(\proback{s})\,\dif s
\end{equation}
\end{lemma}
\begin{proof} The structure of the proof is identical to that of lemma \ref{Xdown_satisfies_SPDE}, so we leave the details to the reader.
\end{proof}
We can hence apply theorem \ref{convergence_of_solution_of_2nd_order_SPDEs} and claim that ${\bf P}[\lim_{t\to\infty}\proback{t}\text{ exists for all }x]=1$. The limit is a random variable measurable in $\forfilt{-\infty}{0}$ independent of $x$ whose distribution is the invariant measure $\invmeas$.
\begin{corollary}
Let \proup{t} a positive recurrent diffusion as in \eqref{recurrent_diffusion}. Then its random attractor is one single point.
\end{corollary}
\begin{proof}
Following \cite{arnold}, the random attractor is defined as the closed random set $\mathcal{A}(\mathfrak{b})\in\forfilt{-\infty}{0}$ such that 
\begin{enumerate}\renewcommand{\theenumi}{\alph{enumi}} 
\item $\flowup{t}(\cdot, \bp)\circ\mathcal{A}(\bp)=\mathcal{A}(\theta_t\bp)$ for all $t\ge0$ almost surely;
\item For all bounded set $D\in\real{}$,
\[
{\bf P}\left[\lim_{t\to\infty}\mathrm{dist}\left(\flowup{t}\big(D, \theta_{-t}\bp),\mathcal{A}(\mathfrak{b})\right)=0\right]=1.
\]
\end{enumerate}
(a) follows from the composition of maps rule. (b) is a straightforward consequence of lemma \ref{spde_for_pull-back_process} and theorem \ref{convergence_of_solution_of_2nd_order_SPDEs}.
\end{proof}
%%%%%%%%%%%%%%%%%%%%%%%%%%%%%%%%%%%%%%%
%% Bibliography
%%%%%%%%%%%%%%%%%%%%%%%%%%%%%%%%%%%%%%%
\bibliographystyle{amsplain}
\bibliography{on_the_invariant_measure_1d}
\end{document}